\documentclass[12pt]{amsart}
\usepackage{amssymb}
\newtheorem{theorem}{Theorem}[section]

\newtheorem{proposition}[theorem]{Proposition}

\newtheorem{lemma}[theorem]{Lemma}
\theoremstyle{definition}

\theoremstyle{remark}
\newtheorem{remark}[theorem]{Remark}
\newtheorem{example}[theorem]{Example}
\newcommand\M{\mathcal{M}}

\newcommand{\V}{\mathcal{V}}

\renewcommand{\O}{\mathcal{O}}

\newcommand{\F}{\mathcal{F}}
\newcommand{\N}{\mathcal{N}}
\newcommand{\R}{\mathbb{R}}
\newcommand{\C}{\mathbb{C}}
\newcommand{\cC}{\mathcal{C}}

\newcommand{\Z}{\mathbb{Z}}
\newcommand{\Q}{\mathbb{Q}}
\renewcommand{\P}{\mathbb{P}}

\newcommand\lie[1]{\mathfrak{#1}}

\renewcommand{\t}{\lie{t}}
\newcommand{\Alc}{\lie{A}}

\newcommand{\su}{\lie{su}}
\newcommand{\on}{\operatorname}

\newcommand{\Ad}{ \on{Ad} }

\newcommand{\Vol}{  \on{Vol}}

\newcommand\dirac{/\kern-1.2ex\partial} 
\newcommand\qu{/\kern-.7ex/} 
\newcommand{\fus}{\circledast}  
\newcommand{\Waff}{W_{\on{aff}}} 

\newcommand{\labell}\label

\newcommand{\ra}{\rightarrow}

\renewcommand{\d}{{\mbox{d}}}
\newcommand{\ol}{\overline}

\newcommand\lam{\lambda}

\newcommand\sig{\sigma}
\newcommand\eps{\epsilon}

\newcommand{\Del}{\Delta}

\renewcommand{\l}{\langle}
\renewcommand{\r}{\rangle}

\newcommand{\ti}{\tilde}

\newcommand\cE{\mathcal{E}}

\newcommand\cF{\mathcal{F}}
\newcommand\cQ{\mathcal{Q}}
\newcommand\pardeg{\on{pardeg}}
\newcommand\rk{\on{rk}}

\renewcommand{\ss}{\on{ss}}
\newcommand\Tr{\on{Tr}}
                                 
\begin{document}

\title[Eigenvalues of products and quantum Schubert calculus] {
Eigenvalues of products of unitary matrices and quantum Schubert
calculus}

\author{S. Agnihotri}
\address{University of Texas at Austin, Department of Mathematics, Austin TX 78712}
\email{sharad@math.utexas.edu}

\author{C. Woodward}
\address{Harvard University, Department of Mathematics,
1 Oxford Street, Cambridge, Massachusetts 02138}
  
\date{\today}

\email{woodward@math.harvard.edu}

\begin{abstract}  We describe the inequalities on the possible eigenvalues
of products of unitary matrices in terms of quantum Schubert calculus.
Related problems are the existence of flat connections on the
punctured two-sphere with prescribed holonomies, and the decomposition
of fusion product of representations of $SU(n)$, in the large level
limit.

In the second part of the paper we investigate how various aspects of
the problem (symmetry, factorization) relate to properties of the
Gromov-Witten invariants.
\end{abstract}

\maketitle


\section{Introduction}

Beginning with Weyl \cite{we:ei}, many mathematicians have been
interested in the following question: given the eigenvalues of two
Hermitian matrices, what are the possible eigenvalues of their sum?
In a recent preprint \cite{kl:sb}, Klyachko observes that a complete
solution to this problem is given by an application of Mumford's
criterion in geometric invariant theory.  The eigenvalue inequalities
are derived from products in Schubert calculus.  In particular, Weyl's
inequalities correspond to Schubert calculus in projective space.  The
necessity of these conditions is due to Helmke and Rosenthal
\cite{hr:ei}.

One of the fascinating points about the above problem are several
equivalent formulations noted by Klyachko.  For instance, the problem
is related to the following question in representation theory: Given a
collection of irreducible representations of $SU(n)$, which
irreducibles appear in the tensor product?  A second equivalent
problem involves toric vector bundles over the complex projective
plane.

In this paper we investigate the corresponding problem for {\em
products} of {\em unitary} matrices.  This question also has a
relationship with a representation-theoretic problem, that of the
decomposition of the fusion product of representations.  The solution
to the multiplicative problem is also derived from geometric invariant
theory, namely from the Mehta-Seshadri theory of parabolic bundles
over the projective line.  The main result of this paper, Theorem
\ref{final}, shows that the eigenvalue inequalities are derived from
products in {\em quantum} Schubert calculus.  This improves a result
of I. Biswas \cite{bi:ex}, who gave the first description of these
inequalities.  A similar result has been obtained independently by
P. Belkale \cite{bl:ip}.

The proof is an application of the Mehta-Seshadri theorem.  A set of
unitary matrices $A_1,\ldots ,A_l$ such that each $A_i$ lies in a
conjugacy class $\cC_i$ and such that their product is the identity is
equivalent to a unitary representation of the fundamental group of the
$l$ times punctured sphere, with each generator $\gamma_i$ being mapped to
the conjugacy class $\cC_i.$ By the Mehta-Seshadri theorem such
a representation exists if and only if there exists a semi-stable
parabolic bundle on $\P^1$ with $l$ parabolic points whose parabolic
weights come from the choice of conjugacy classes $\cC_i.$  This last
interpretation of the original eigenvalue problem can be related to
the Gromov-Witten invariants of the Grassmannian and this is done in
Section 5 below.

In Sections 6 and 7 we investigate how factorization and hidden
symmetries of these Gromov-Witten invariants relate to the
multiplicative eigenvalue problem.

\section{Additive inequalities (after Klyachko and Helmke-Rosenthal)}

Let $\su(n)$ denote the Lie algebra of $SU(n)$, and 
$$\t = \{ (\lambda_1,\ldots,\lambda_n) \in \R^n \ | \
\sum \lambda_i = 0 \} $$
its Cartan subalgebra.  Let 
$$ \t_+ =  \{ (\lambda_1,\ldots,\lambda_n) \in \t \ | \ 
 \lam_i  \ge \lam_{i+1}, \ \ i=1,\ldots,n-1 \} $$
be a choice of closed positive Weyl chamber.  For any matrix $A \in
\su(n)$ let
$$ \lam(A) = (\lam_1(A),\lam_2(A),\cdots,\lam_n(A)) \in \t_+ $$
be the eigenvalues of the Hermitian matrix $-i A$ in non-increasing
order.  Let $\Del(l) \subset (\t_+)^l$ denote the set
$$ \Del(l) = \{ (\lam(A_1),\lam(A_2),\ldots,\lam(A_l)) \ \vert \ 
A_1,\ldots,A_l \in \su(n), \ A_1 + A_2 + \ldots + A_l = 0 \} .$$
Define an involution
$$ *: \ \t_+ \cong \t_+, \ \ (\lambda_1,\ldots,\lambda_n)
\mapsto (-\lambda_n, \ldots, -\lambda_1) .$$
For any $A \in \su(n)$ the matrix $-A$ has eigenvalues $ \lam(-A) =
*\lam(A)$.  The set $\Del(l)$ is invariant under the map
$$ *^l:\  (\t_+)^l \to (\t_+)^l, \ \ (\xi_1,\ldots,\xi_l) \mapsto
(*\xi_1,\ldots,*\xi_l) $$
and also under the action of the symmetric group $S_l$ on
$(\t_+)^l$.

The set $\Del(l)$ has interesting interpretations in symplectic
geometry and representation theory.  Consider the cotangent bundle $
T^*SU(n)^{l-1} $ with the action of $SU(n)^l$ given by $SU(n)$ acting
diagonally on the left and $SU(n)^{l-1}$ on the right.  The moment
polytope of this action may be identified with $\Del(l)$ (see Section
5.)  From convexity theorems in symplectic geometry (see e.g.
\cite{sj:co} and \cite{le:co}) it follows that $\Del(l)$ is a
finitely-generated convex polyhedral cone.  In particular there are a
finite number of inequalities defining $\Del(l) $ as a subset of the
polyhedral cone $(\t_+)^l$.

The set $\Del(l) $ may also be described in terms of the tensor
product of representations. Let 
$$(\ ,\ ) : \ \su(n) \times \su(n) \to \R, \ \ \ (A,B) \mapsto - \Tr(AB) $$ 
denote the basic inner product on $\su(n)$, which induces an
identification $\su(n) \cong \su(n)^*$.  Let $ \Lambda = \Z^n \cap \t$
denote the integral lattice and $ \Lambda^* \subset \t$ its dual, the
weight lattice.  For each $\lam \in \Lambda^*\cap \t_+$, let $V_\lam$
denote the corresponding irreducible representation of $SU(n)$.  We
will see in equation \eqref{inv_eqn} that $\Del(l) \cap \Q^l$ is the
set of $(\lam^1,\ldots,\lam^l)$ such that for some $N$ such that
$N\lam^1,\ldots, N \lam^l \in \Lambda^*$, we have
$$V_{N \lam^1} \otimes \ldots \otimes V_{N\lam^{l-1}} \supset V^*_{N
\lam^l},$$
that is, $V_{N \lam^1} \otimes \ldots \otimes V_{N \lam^l}$ contains
a non-zero invariant vector.

The work of Klyachko and Helmke-Rosenthal gives a complete set of
inequalities describing $\Del(l)$ in terms of Schubert calculus.  Let
$$ \C^n = F_n \supset F_{n-1} \supset \ldots \supset F_0 = \{ 0 \} $$ 
be a complete flag in $\C^n$, $G(r,n)$ the Grassmanian of $r$-planes
in $\C^n$, and for any subset $I = \{i_1,\ldots,i_r \} \subset \{1,
\ldots n \}$ let
$$ \sig_I = \{ W \in G(r,n) \ \vert \ \dim(W \cap F_{i_j}) \ge j, \ 
j = 1,\ldots,r \} $$
denote the corresponding Schubert variety.  The Schubert cell $C_I
\subset \sigma_I$ is defined as the complement of all lower-dimensional Schubert varieties contained in $\sig_I$:
$$ C_I = \sigma_I \backslash \bigcup_{\sigma_J \subset \sigma_I}
\sigma_J.
$$
We say that $W$ is in {\em position $I$} with respect to the flag
$F_*$ if $W \in C_I$.

The homology classes $[ \sig_I ]$ form a basis of $H_*(G(r,n),\Z)$.
Given two Schubert cycles $\sig_I,\sig_J$, we can expand the
intersection product $[\sig_I] \cap [\sig_J]$ in terms of this basis.
We say $[\sig_I] \cap [\sig_J]$ contains $[\sig_K]$ if $[\sig_K]$
appears in this expansion with non-zero (and therefore positive)
coefficient.  Equivalently, let
$$*K = \{n + 1 - i_r,n+ 1 - i_{r-1},\ldots, n+1-i_1 \},$$
so that $[\sig_{*K}]$ is the Poincare dual of $[\sig_K]$.  Then
$[\sig_I] \cap [\sig_J]$ contains $[\sig_K]$ if and only if the
intersection of general translates of the Schubert cycles
$\sig_I,\sig_J,\sig_{*K}$ is non-empty and finite.

\begin{theorem}[Klyachko, resp. Helmke-Rosenthal]  \labell{klyachko}
  A complete (resp. necessary) set of inequalities
  describing $\Del(l)$ as a subset of $(\t_+)^l$ are
\begin{equation} \labell{add_ineq}
\sum_{i \in I_1} \lam_i(A_1) + \sum_{i \in I_2} \lam_i(A_2)
+ \ldots + \sum_{i \in I_l} \lam_i(A_l) \leq 0,
\end{equation}
where $I_1,\ldots,I_l$ are subsets of $\{ 1 ,\ldots, n \}$ of the same
cardinality $r$ such that $[\sig_{I_1}] \cap \ldots \cap
[\sig_{I_{l-1}}] \supset [\sig_{*I_l}]$, and $r$ ranges over all
values between $1$ and $n-1$.
\end{theorem}

Note that the cases $l =1,2$ are trivial: $\Del(1) = \{ 0 \},$ and $
\Del(1) = \{ (\mu,*\mu) \ |\ \mu \in \t_+ \}.$ Klyachko also claims
that these inequalities are independent.  From Theorem \ref{klyachko}
follows a complete set of inequalities for the possible eigenvalues of
a sum of skew-Hermitian matrices.  For instance, for $l=3$ one obtains
the inequalities
\begin{equation} \labell{add_ineq2}
\sum_{i \in I}
\lam_i(A) + \sum_{j \in J} \lam_j(B) \le \sum_{k \in K} \lam_k(A + B),
\end{equation}
where $I,J,K \subset \{ 1 ,\ldots, n \}$ range over subsets such that
$[\sig_I] \cap [\sig_J]$ contains $[\sig_K]$.

\begin{example}  Let $r=1$ so that $G(r,n) \cong \P^{n-1}$ and
$I = \{ n - i + 1 \}, \ J = \{ n - j + 1 \}$.  Then $\sigma_I \cong
\P^{n-i},\ \sigma_J \cong \P^{n-j}$ so that $[\sigma_I] \cap
[\sigma_J] \cong \P^{n-i-j+1} = \sigma_K$ where $K = \{n - i -j +
2\}$.  One obtains
\begin{equation} \labell{dual_Weyl}
\lam_{n-i+1}(A) + \lam_{n-j+1}(B) \le \lam_{n-i-j+2}(A + B) .
\end{equation}
\end{example}

\subsection{Duality}  Let $A_1,\ldots, A_l \in \su(n)$.    
From \eqref{add_ineq2} applied to $-A_1,\ldots ,- A_l$ one obtains
\begin{equation} \labell{dual_ineq}
 - \sum_{i \in *I_{1}} \lam_i(A_1) -\ldots - \sum_{i \in *I_{l}} \lam_i(A_l)\leq 0
 \end{equation} 
or equivalently $   \sum_{i \in *I_{1}} \lam_i(A_1) +\ldots + \sum_{i \in *I_{l}} \lam_i(A_l)\geq 0. $
By the trace condition, \eqref{dual_ineq} is equivalent to
$$    \sum_{i \notin *I_{1}} \lam_i(A_1) +\ldots + \sum_{i \notin *I_{l}} \lam_i(A_l)\leq 0  .$$
Let $I_i^c = \{ 1,\ldots,n \} \backslash * I_i$.  Then
$[\sigma_{I_i^c}]$ is the image of $[\sigma_{I_{i}}]$ under the
isomorphism of homology induced by $G(r,n) \cong G(n-r,n)$ (see page
197 onwards of Griffiths and Harris \cite{gr:pr}).  Thus the
appearance of \eqref{dual_ineq} in \eqref{add_ineq} corresponds to a
product in the Schubert calculus of $G(n-r,n)$.

\begin{example} \labell{Weyl_example}
The dual equation to \eqref{dual_Weyl} is Weyl's 1912 \cite{we:ei}
inequality
\begin{equation} \labell{Weyl} 
\lam_i(A) + \lam_j(B) \ge \lam_{i + j -1}(A + B).
\end{equation}
\end{example}

\section{Multiplicative Inequalities}

Let $\Alc \subset \t_+$ be the fundamental alcove of $SU(n)$: 
$$\Alc = \{ \lam \in \t_+ \ | \ \lam_1 - \lam_n \leq 1 \}.$$
Let $A \in SU(n)$ be a unitary matrix with determinant $1$.  Its
eigenvalues may be written
$$ e^{2\pi i \lam_1(A)}, e^{2\pi i \lam_2(A)}, \ldots ,e^{2\pi i \lam_n(A)} $$
where $\lam(A) = (\lam_1(A),\ldots,\lam_n(A)) \in \Alc$.  The map $A
\mapsto \lam(A)$ induces a homeomorphism
$$ \Alc \cong SU(n) / \Ad(SU(n)) .$$

Let $\Del_q(l) \subset \Alc^l$ ($q$ for quantum) denote the set
$$ \Del_q(l) = \{ (\lam(A_1),\ldots,\lam(A_l))  \ \vert \ 
A_1,\ldots,A_l \in SU(n), \ A_1A_2\ldots A_l = I \}. $$
As before, $\Del_q(l)$ is invariant under the involution, $*^l: \Alc^l
\to \Alc^l$, and the action of the symmetric group $S_l$ on $\Alc^l$.

The set $\Del_q(l)$ has an interpretation as a moment polytope.  Let
$\M$ be the space of flat $SU(n)$-connections on the trivial $SU(n)$
bundle over the $l$-holed two-sphere, modulo gauge transformations
which are the identity on the boundary (see \cite{me:lo}).  The gauge
group of the boundary acts on $\M$ in Hamiltonian fashion and the set
$\Del_q(l)$ is the moment polytope for this action.  By \cite[Theorem
3.19]{me:lo}, $\Del_q(l)$ is a convex polytope. In fact, an analogous
statement holds for arbitrary compact, simply-connected Lie groups.
In particular, a finite number of inequalities describe $\Del_q(l)$.
In the case $n=2$, these inequalities were given explicitly for $l=3$
in Jeffrey-Weitsman \cite{jw:bs} and for arbitrary numbers of marked
points in Biswas \cite{bi:r2}.  A description of the inequalities in
the arbitrary rank case was given in \cite{bi:ex} but the description
given there does not seem to be computable.

There is also an interpretation of $\Del_q(l)$ in terms of fusion
product.  Let $\fus_N$ denote the fusion product on the Verlinde
algebra $R(SU(n)_N)$ of $SU(n)$ at level $N$.  Then $\Del_q(l) \cap
\Q^l$ is the set of $(\lam^1,\ldots,\lam^l) \in \Alc \cap \Q^l$ such
that for some $N$ such that $N\lam^1,\ldots, N \lam^l \in \Lambda^*$,
we have 
\begin{equation} \labell{fus}
V_{N \lam^1} \fus_N \ldots \fus_N V_{N\lam^{l-1}} \supset
V_{N * \lam^l}.\end{equation}
See Section \ref{Verlinde}.

\subsection{Quantum Schubert calculus}

Quantum cohomology is a deformation of the ordinary cohomology ring
that was introduced by the physicists Vafa and Witten.  Quantum
cohomology of the Grassmannian (quantum Schubert calculus) was put on
a rigorous footing by Bertram \cite{be:qs}.  Recall that the degree of
a holomorphic map $\varphi: \ \P^1 \ra G(r,n)$ is the homology class
$[\varphi] \in H^2(G(r,n),\Z) \cong \Z$.  Let $p_1,\ldots,p_l$ be
distinct marked points in $\P^1$.  The quantum intersection product
$\star$ on $H_*(G(r,n),\C) \otimes \C[q]$ is defined by
$$ [\sigma_{I_1} ] \star \ldots \star [\sigma_{I_l}] = \sum_{J}
\langle[\sig_{I_1}],\ldots, [\sigma_{I_l}], [\sig_J] \rangle_d \
[\sigma_{*J}] q^d,$$
where the Gromov-Witten invariant $\langle[\sig_{I_1}],\ldots,
[\sigma_{I_l}], [\sig_J] \rangle_d$ is equal to the number of
holomorphic maps $\P^1 \mapsto G(r,n)$ sending $p_1,\ldots,p_l,p$ to
general translates of $\sigma_{I_1},\ldots,\sigma_{I_l},\sigma_J$ if
this number is finite, and is otherwise zero.

Our main result is the following description of $\Del_q(l)$:
\begin{theorem}\labell{final}   A complete set of inequalities for 
  $\Del_q(l)$ are given by 
\begin{equation}
\sum_{i \in I_1} \lam_i(A_1) + \sum_{i \in I_2} \lam_i(A_2)
+ \ldots + \sum_{i \in I_l} \lam_i(A_l) \leq d
\end{equation}
for $(I_1,\ldots,I_l,d)$ such that $\l [\sig_{I_1}] \star \ldots \star
[ \sig_{I_l} ] \r_d \neq 0,$ that is, $\ [\sig_{I_1}] \star \ldots
\star [ \sig_{I_{l-1}}] \supset q^d [\sig_{*I_l}] .$
\end{theorem}

In the last few years several techniques have been developed for
computing the coefficients of quantum Schubert calculus.  See for
instance Bertram, Ciocan-Fontanine, Fulton \cite{be:qm}.  Therefore
the above theorem makes the question of which inequalities occur
computable in practice.

One recovers the inequalities for $\Del(l)$ from the degree $0$
Gromov-Witten invariants.  This shows that $\Del(l) $ is
the cone on $\Del_q(l)$ at the $0$-vertex, i.e.
$$ \Del(l) = \R_+ \cdot \Del_q(l) $$ 
This may be verified by several alternative methods, e.g. 
Remark \ref{small}.

The simplest example of a positive degree inequality is given by the
following:
\begin{example} \labell{quantum_Weyl_example}  
Let $r = 1$ so that $G(r,n) = \P^{n-1}$, and $U,V,W \subset \C^n$ be
subspaces in general position of dimensions $i,j,n+1-i-j$.  There is a
unique degree $1$ map $\P^1 \ra \P^{n-1}$ mapping $p_1,p_2,p_3$ to
$\P(U),\P(V),\P(W)$ respectively.  Together with the degree $0$
inequality mentioned before, this gives
\begin{equation} \labell{quant_Weyl}
 \lam_{i+j-1}(AB) \le \lam_i(A) + \lam_j(B) \le \lam_{i+j}(AB) + 1.
\end{equation}
We will see in Section \ref{symmetry} that these inequalities are
related by a symmetry of $\Del_q(l)$.
\end{example}

\vskip .1in {\noindent \em Question:} Are the inequalities in Theorem
\ref{final} independent?  \vskip .1in

\section{Moduli of flags and Mumford's criterion}

As a warm-up we review some of the ideas involved in Klyachko's
proof.  For any $\xi \in \t_+$ let
$$\O_\xi = SU(n) \cdot \xi = \{ A \in \su(n) \ | \ \lam(A) = \xi \} $$
denote the corresponding adjoint orbit.  Via the identification
$\su(n) \cong \su(n)^*$, $\O_\xi$ inherits a canonical symplectic
structure, called the Kostant-Kirillov-Souriau two-form, and the
action of $SU(n)$ on $\O_\xi$ is Hamiltonian with moment map given by
inclusion into $\su(n)$.

The diagonal action of $SU(n)$ on $\O_{\xi_1} \times \ldots \times
\O_{\xi_l}$ has moment map given by
$$ (A_1,\ldots,A_l) \mapsto A_1 + A_2 + \ldots + A_l.$$
The symplectic quotient $\N({\xi_1,\ldots,\xi_l}) = \O_{\xi_1} \times
\ldots \times \O_{\xi_l} \qu SU(n)$ is given by
$$ \N({\xi_1,\ldots,\xi_l}) = \{(A_1,\ldots,A_l) \in \O_{\xi_1} \ldots
\O_{\xi_l} \ | \ A_1 + A_2 + \ldots + A_l = 0 \} / SU(n) .$$
For generic $(\xi_1,\ldots,\xi_l),$ that is, values where the moment
map has maximal rank, the quotient $\N({\xi_1,\ldots,\xi_l})$ is a
symplectic manifold.  The $l$-tuple $(\xi_1,\ldots,\xi_l)$ lies in
$\Del(l) $ if and only if $\N({\xi_1,\ldots,\xi_l})$ is non-empty.

The quotients $\N({\xi_1,\ldots,\xi_l})$ may be viewed as symplectic
quotients of the cotangent bundle $ T^*SU(n)^{l-1} $.  Indeed, the
symplectic quotient 
$$ (T^*SU(n) \times \O_\xi) \qu SU(n) \cong \O_\xi .$$
Therefore, the quotient of $T^*SU(n)^{l-1}$ by the right action of
  $SU(n)^{l-1}$ and the diagonal left action of $SU(n)$ is
  $$ (T^*SU(n)^{l-1} \times \O_{\xi_1} \times \ldots \times \O_{\xi_{l}})
  \qu SU(n)^{l} \cong \N({\xi_1,\ldots,\xi_l}).$$
It follows that $\Del(l)$ is the moment polytope of the action of
$SU(n)^l$ on $ T^*SU(n)^{l-1}$.

One can determine whether $\N({\xi_1,\ldots,\xi_l})$ is empty by
computing its symplectic volume.  This is given by a formula derived
from the Duistermaat-Heckman theorem due to Guillemin-Prato (see
\cite{gu:he} or \cite[(4)]{me:wi}).  Unfortunately the formula involves
cancelations and it is not apparent what the support of the volume
function is, or even that the support is a convex polytope.

The manifolds $\O_{\xi_i}$ have canonical complex structures (induced
by the choice of positive Weyl chamber) and are isomorphic to
(possibly partial) flag varieties.  Suppose that $\xi_1,\ldots,\xi_l$
lie in the weight lattice $\Lambda^*$, so that there exist pre-quantum
line bundles $L_{\xi_i} \to \O_{\xi_i}$; i.e., equivariant line
bundles with curvature equal to $2 \pi i$ times the symplectic
form. The sections of $L_{\xi_1}\boxtimes \ldots \boxtimes L_{\xi_l}$
define a K\"ahler embedding 
$$\O_{\xi_1} \times \ldots \times
\O_{\xi_l} \to \P(V_{\xi_1}^*) \times \ldots \times \P(V_{\xi_l}^*),$$
where $V_{\xi_1},\ldots,V_{\xi_l}$ are the irreducible representations
with highest weights $\xi_1,\ldots,\xi_l$.  By an application of a
theorem of Kirwan and Kempf-Ness (which holds for arbitrary smooth
projective varieties, see \cite[page 109]{ki:coh}) the symplectic
quotient is homeomorphic to the geometric invariant theory quotient
$$\N({\xi_1,\ldots,\xi_l}) \cong \O_{\xi_1} \times \ldots \times \O_{\xi_l} \qu
SL(n,\C).$$
By definition, $\O_{\xi_1} \times \ldots \times \O_{\xi_l} \qu
SL(n,\C) $ is the quotient of the set of semi-stable points in
$\O_{\xi_1} \times \ldots \times \O_{\xi_l}$ by the action of
$SL(n,\C)$, where $(F_{\xi_1},\ldots,F_{\xi_l}) \in \O_{\xi_1} \times
\ldots \times \O_{\xi_l}$ is called semi-stable if and only if for
some $N$ there is an invariant section of $(L_{\xi_1} \boxtimes \ldots
\boxtimes L_{\xi_l})^{\otimes N} $ which is non-vanishing at
$(F_{\xi_1},\ldots,F_{\xi_l})$.  The quotient $\N({\xi_1,\ldots,\xi_l})$
is therefore non-empty if and only if there exists a non-zero
$SU(n)$-invariant vector in
\begin{equation} \labell{inv_eqn}
 H^0((L_{\xi_1} \boxtimes \ldots \boxtimes L_{\xi_l})^{\otimes N}) =
V_{N \xi_1} \otimes \ldots \otimes V_{N \xi_l}.\end{equation}
This explains the representation-theoretic interpretation of $\Del(l)
$ alluded to in the introduction.

In order to obtain the inequalities in Theorem \ref{klyachko}, one
applies the criterion of Mumford, which says that {\em a point is
semi-stable if and only if it is semi-stable for all one-parameter
subgroups} \cite[Chapter 2]{mu:ge}, see also \cite[Lemma 8.8]{ki:coh}.
Let us assume that the $\xi_i$ are generic.  An application of the
criterion gives that an $l$-tuple of complete flags
$(F_{1},\ldots,F_{l}) \in \O_{\xi_1} \times \ldots \times \O_{\xi_l}$
is semi-stable if and only if for all subspaces $W \subset \C^n$, one
has
$$ \sum_{i \in I_1} \xi_{1,i} + \ldots + \sum_{i \in I_l} \xi_{l,i} \le 0, $$
where $I_j$ is the position of $W$ with respect to the flag $F_{j}$.
The proof similar to that for Grassmannians given in Section 4.4 of
\cite{mu:ge}.

The set of semi-stable points is dense if non-empty.  It follows that
$\N({\xi_1,\ldots,\xi_l})$ is non-empty if and only if the above
inequality holds for every intersection $\sigma_{I_1} \cap \ldots \cap
\sigma_{I_l}$ of Schubert cycles in general position.  Any inequality
corresponding to a positive dimensional intersection must be
redundant.  Indeed, since the intersection is a projective variety, it
cannot be contained in any of the Schubert cells.  The boundary of
$\sigma_{I_l}$ consists of Schubert varieties $\sigma_J$ with $J$ such
that $j_k \le i_k$ for $k = 1,\ldots,r$, where $i_1,\ldots,i_r$ and
$j_1,\ldots, j_r$ are the elements of $I_l$ and $J$ in increasing
order.  The inequality obtained from an intersection $\sigma_{I_1}
\cap \ldots \cap \sigma_{I_{l-1}} \cap \sigma_J \neq \emptyset$
therefore implies the inequality obtained from $\sigma_{I_1} \cap
\ldots \cap \sigma_{I_l} \neq \emptyset$.

\section{Application of the Mehta-Seshadri theorem}

For any $\xi \in \Alc$, let
$$\cC_{\xi} = \{ A \in SU(n) \ | \ \lam(A) = \xi \} $$
denote the corresponding conjugacy class.  The mapping $A \mapsto
\lam(A)$ induces a homeomorphism $SU(n)/\Ad(SU(n)) \cong \Alc $.

Let $p_1,\ldots,p_l \in \P^1$ be distinct marked points and
$\M(\xi_1,\ldots,\xi_l)$ the moduli space of flat $SU(n)$-connections
on $\P^1 \backslash \{p_1,\ldots, p_l \} $ with holonomy around $p_i$
lying in $\cC_{\xi_i}$.  Since the fundamental group of $\P^1
\backslash \{p_1,\ldots, p_l\}$ has generators the loops $\gamma_1,
\ldots, \gamma_l$ around the punctures, with the single relation
$\gamma_1 \cdot \ldots \cdot \gamma_l=1$,
$$ \M({\xi_1,\ldots,\xi_l}) \cong \{ (A_1,\ldots,A_l) \in \cC_{\xi_1}
\times \ldots \times \cC_{\xi_l} \ | \ A_1A_2 \cdots A_l = I \} /
SU(n) .$$
In particular $\M({\xi_1,\ldots,\xi_l})$ is non-empty if and only if
$(\xi_1,\ldots,\xi_l) \in \Del_q(l)$. In theory one can determine if
$\M({\xi_1,\ldots,\xi_l})$ is non-empty by computing its symplectic
volume by the formulae stated in Witten \cite[(4.11)]{wi:tw}, Szenes
\cite{sz:vo}, and \cite[Theorem 5.2]{me:co}.

For rational $\xi_1,\ldots,\xi_l$ the space $\M({\xi_1,\ldots,\xi_l})$
has an algebro-geometric description due to Mehta-Seshadri
\cite{ms:pb}.  Let $C$ be a Riemann surface with marked points
$p_1,\ldots,p_l \in C$ and let $\cE \ra C$ be a holomorphic bundle.  A
parabolic structure without multiplicity on $\cE$ consists of the
following data at each marked point $p_i$: a complete ascending flag
$$ 0 = \cE_{p_i,0} \subset \cE_{p_i,1} \subset \cE_{p_i,2} \ldots \subset \cE_{p_i,n} = \cE_{p_i} $$
in the fiber $\cE_{p_i}$ and a set of {\em parabolic weights}
$$ \lam_{i,1} > \lam_{i,2} > \ldots > \lam_{i,n} $$
satisfying $\lam_{i,1} - \lam_{i,n} \leq 1$.  In \cite{ms:pb} the
weights are required to lie in the interval $[0,1)$, but the
definitions work without this assumption.  A parabolic bundle is a
holomorphic bundle with a parabolic structure.  Recall that the degree
$\deg(\cE)$ of $\cE$ is the first Chern class $c_1(\cE) \in H^2(C,\Z)
\cong \Z$.  The parabolic degree $\pardeg(\cE)$ is defined by
$$ \pardeg(\cE) = \deg(\cE) + \sum_{i=1,j=1}^{l,n} \lam_{i,j} .$$
The parabolic slope $\mu(\cE)$ is
$$ \mu(\cE) = \frac{\pardeg(\cE)}{\rk(\cE)} .$$
Given a holomorphic sub-bundle $\F \subset \cE$ of rank $r$ one
obtains a parabolic structure on $\F$ as follows.  An ascending flag
in the fiber $\F_{p_i}$ at each marked point $p_i$ is obtained by
removing from
$$ \F_{p_i} \cap \cE_{p_i,1} \subseteq \F_{p_i} \cap \cE_{p_i,2}
 \subseteq \ldots \subseteq \F_{p_i} \cap \cE_{p_i,n} = \F_{p_i} $$
those terms for which the inclusion is not strict.  The parabolic
weights for $\F$ are
$$ \mu_{i,j} = \lam_{i,k_j}, $$
where $k_j$ is the minimal index such that $\F_{p_i,j} \subseteq
\cE_{p_i,k_j}$.  Let $K_i = \{ k_1, \ldots, k_r \}$.  The fiber
$\F_{p_i}$ may be viewed as a element of the Grassmannian of
$r$-planes in $\cE_{p_i}$, and $K$ is the position of $\F_{p_i}$ with
respect to the flag $\cE_{p_i,*}$.  The parabolic degree of $\F$ is
$$ \pardeg(\F) = \deg(\F) + \sum_{i,\ k \in K_i} \lam_{i,k} .$$

A parabolic sub-bundle of $\cE$ is a holomorphic sub-bundle $\F
\subset E$ whose parabolic structure is the one induced from the
inclusion.  A parabolic bundle $\cE \to C$ is called parabolic
semi-stable if $\mu(\F) \le \mu(\cE)$ for all parabolic sub-bundles
$\F \subset \cE$.  There is a natural equivalence relation on
parabolic bundles: Two bundles are said to be grade equivalent if the
associated graded bundles are isomorphic as parabolic bundles.  See
\cite{ms:pb} for more details.

\begin{theorem}[Mehta-Seshadri]   Suppose the parabolic weights $\lam_{i,j}$
are rational and lie in the interval $[0,1)$.  Then the moduli space
$\M(\lam_1, \ldots,\lam_l)$ of grade equivalence classes of
semi-stable parabolic bundles with parabolic weights $\lam_{i,j}$ and
parabolic degree $0$ is a normal, projective variety, homeomorphic to
the moduli space of flat unitary connections over $C \backslash \{
p_1,\ldots,p_r \}$ such that the holonomy of a small loop around $p_i$
lies in $\cC_{\lam_{i}}$.
\end{theorem}

In fact, the Mehta-Seshadri theorem also holds without the assumption
that the parabolic weights lie in $[0,1)$.  One can see this either
through the theory of elementary transformations, or through the
extension of the Mehta-Seshadri theorem to non-zero parabolic degree
given in Boden \cite{bo:re}.

The explanation using elementary transformations goes as follows. Let
$\cQ$ denote the skyscraper sheaf with fiber $\cE_{p_i}/\cE_{p_i,n-1}$
at $p_i$.  One has an exact sequence of sheaves
$$ 0 \to \cE' \to \cE \to \cQ \to 0 .$$
The kernel $\cE'$ is a sub-sheaf of a locally free sheaf and therefore
locally free.  Since degree is additive in short exact sequences
$\deg(\cE') = \deg(\cE)-1$.  One calls the $\cE'$ an elementary
transformation of $\cE$ at $p_i$.  There is a canonical line
$\cE'_{p_i,1}$ in the fiber $\cE'_{p_i}$ which is the kernel of the
fiber map $\pi: \cE'_{p_i} \to \cE_{p_i}$.  One extends the canonical
line to a complete flag by taking $\cE'_{p_i,j} =
\pi^{-1}(\cE'_{p_i,j-1})$ for $j>1$.  Finally one takes as parabolic
weights at $p_i$ the set $ \lam_{i,n} + 1,\lambda_{i,1},
\ldots,\lam_{i,n-1}$.  With this parabolic structure the bundle $\cE'$
is parabolic semi-stable of the same parabolic degree as $\cE$.  Details,
in a slightly different form, can be found in Boden and Yokogawa
\cite{bo:ra}.

The following is the key lemma in the derivation of Theorem
\ref{result} from Mehta-Seshadri.  Let $d = \deg(\cE) = -\sum
\lam_{i,j}$ denote the degree of any element $\cE \in \M(\lam_1,
\ldots,\lam_l)$.

\begin{lemma} \labell{key}  Suppose that there is some ordinary semi-stable
bundle on $C$ of degree $d$.  Then the set of equivalence classes of
parabolic semi-stable bundles of parabolic degree $0$ whose underlying
holomorphic bundle is ordinary semi-stable is Zariski dense in
$\M(\lam_1,\ldots,\lam_l)$.
\end{lemma}

\begin{proof}  
Recall from the construction of $\M(\lam_1,\ldots,\lam_l)$ in \cite{ms:pb} that for
some integer $N$ there exists an $SL(N)$-equivariant bundle
$ \ti{R} \stackrel{\pi}{\to} R $
whose fibers are products of  $l$ complete flag varieties, such that the
geometric invariant theory quotients of $\ti{R},R$ are
$$ \ti{R} \qu SL(N) = \M(\lam_1,\ldots,\lam_l) , \ \ \ \ R \qu SL(N) = \M,$$
where $\M$ denotes the moduli space of ordinary semi-stable bundles on
$C$ of degree $d$.  By definition,
$$ \ti{R} \qu SL(N) = \ti{R}^{\ss} / SL(N), \ \ \ \ R \qu SL(N) =
R^{\ss} / SL(N) $$
where $\ti{R}^{\ss},R^{\ss}$ denote the Zariski dense set of
semi-stable points in $\ti{R},R$ respectively.  The inverse
image $\pi^{-1}(R^{\ss}) \cap \ti{R}^{\ss} / SL(N)$ is therefore dense
in $\ti{R}^{\ss} / SL(N) = \M(\lam_1,\ldots,\lam_l)$.
\end{proof}

Now we specialize to the case $C = \P^1$ with $l$ marked points
$p_1,p_2,\ldots,p_l$.  Let $\xi_1,\ldots,\xi_l \in \Alc^l \cap \Q^l$.
By Lemma \ref{key}, $\M({\xi_1,\ldots,\xi_l})$ is non-empty if and
only there exists a parabolic semi-stable $\cE$ with parabolic degree $0$
and weights $\xi_1,\ldots,\xi_l$ whose underlying holomorphic bundle
is semi-stable.  Since the sum of the parabolic weights is zero, the
degree of $\cE$ is also zero.  By Grothendieck's theorem, $\cE$ is
holomorphically trivial.  A sub-bundle $\cF \subset \cE$ of rank $r$
is given by a holomorphic map
$$ \varphi_{\cF}: \ \P^1 \to G(r,n).$$
Since $\varphi_{\cF}$ is the classifying map of the quotient
$\cE/\cF$, the degree of $\cF$ is minus the degree of $\varphi_{\cF}$.
The parabolic slope of $\cF$ is given by
$$ \mu(\cF) = -\deg(\varphi_{\cF}) + \sum_{i \in I_1(\varphi)} \xi_{1,i} +
\ldots +   \sum_{i \in I_1(\varphi)} \xi_{l,i},  $$
where $I_i(\varphi)$ is the position of the subspace $\varphi(p_i)
\subset \cE_{p_i}$ with respect to the flag $\cE_{p_i,* }$ above.
The parabolic bundle $\cE$ is called parabolic semi-stable if and only
if for all such $F$, $ \mu(F) \leq 0 $, that is,
$$ \sum_{i \in I_1(\varphi)} \xi_{1,i} + \ldots + \sum_{i \in
I_1(\varphi)} \xi_{l,i} \leq \deg(\varphi) $$
for all maps $\varphi:\P^1 \to G(r,n)$.  

The following result was obtained independently by P. Belkale \cite{bl:ip}.

\begin{theorem} \labell{result}  
A complete set of inequalities for $\Del_q(l)$ as a subset of $\Alc^l$
is given by
\begin{equation} \labell{mult_ineq}
\sum_{i \in I_1} \lam_i(A_1) + \sum_{i \in I_2} \lam_i(A_2)
+ \ldots + \sum_{i \in I_l} \lam_i(A_l) \leq d
\end{equation}
for subsets $I_1,\ldots,I_l \subset \{ 1 ,\ldots, n \}$ of the same
cardinality $r$ and non-negative integers $d$ such that there exists a
rational map $\P^1 \to G(r,n)$ of degree $d$ mapping $p_1,\ldots,p_l$
to the Schubert cells $C_{I_1},\ldots,C_{I_l}$ in general position.
\end{theorem}

\begin{proof}  
If $\M(\xi_1,\ldots,\xi_l)$ is non-empty, then a trivial bundle with a
general choice of flags will be parabolic semi-stable.  Indeed, by the
above discussion the fiber $\on{Flag}^l$ of $\pi :\ti{R} \to R$ over a
trivial bundle intersects $\ti{R}^{ss}$, so $\ti{R}^{ss}\cap
\on{Flag}^{l}$ is open in $\on{Flag}^{l}$.  Therefore,
$\M(\xi_1,\ldots,\xi_l)$ is non-empty if and only if
$$  \sum_{i \in I_1} \xi_{1,i} + \ldots + \sum_{i \in
I_l} \xi_{l,i} \leq d $$
for all subsets $I_1,\ldots,I_l$ and integers $d$ such that there
exists a degree $d$ map sending $p_1,\ldots,p_l$ to general translates
of the Schubert cells $C_{I_1},\ldots,C_{I_l}$.  
\end{proof}

\begin{remark}  \labell{small} For sufficiently small parabolic weights
  $\lam_{i,j}$ any parabolic semi-stable bundle on $\P^1$ is
  necessarily ordinary semi-stable of degree $0$, and therefore
  trivial.  It follows that the moduli spaces
  $\M(\lam_1,\ldots,\lam_l)$ and $\N(\lam_1,\ldots,\lam_l)$ are
  isomorphic.  This shows that Klyachko's result is implied by Theorem
  \ref{result}.
\end{remark}

We now show that the existence of the maps described in Theorem
\ref{result} may be detected by Gromov-Witten invariants.
Let $\sig_{I_1},\ldots,\sig_{I_l}$ be some collection of Schubert
varieties, and consider the expansion
$$ [\sig_{I_1}] \star [\sig_{I_2}] \ldots \star [\sig_{I_l}] =
\sum_i q^i \alpha_i $$
where $\alpha_i \in H_*(G(r,n))$.  (Question: is this product always
non-zero?)  We say that $q^d $ divides $ [\sig_{I_1}] \star
[\sig_{I_2}] \ldots \star [\sig_{I_l}]$ if $\alpha_i = 0$ for all $i <
d$.  The following lemma is stated in Ravi \cite{ra:in}.
\begin{lemma} \labell{compute}  Let $d$ be the lowest degree of a map $\P^1
  \to G(r,n)$ sending $p_1,\ldots,p_l$ to general translates of
  $\sig_{I_1},\ldots,\sig_{I_l}$ respectively.  Then $q^d$ is the
  maximal power of $q$ dividing $[\sig_{I_1}] \star
  \ldots \star [\sig_{I_l}]$.
\end{lemma}

\begin{proof} 
Let $\M_d$ denote the space of maps $\P^1 \to G(r,n)$ of degree $d$,
$\on{ev}^l: \, \M_d \to G(r,n)^l$ the evaluation map, and
$\sig_{I_*}(p_*) = (\on{ev}^l)^{-1}( \sigma_{I_*}) $ the subset of
maps sending $p_j$ to $\sigma_{I_j}$ for $j=1,\ldots,l$.  By
\cite[Moving Lemma 2.2A]{be:qs}, $ \sig_{I_*}(p_*) $ is a
quasi-projective variety, of the expected codimension in $\M_d$.  By
choosing enough additional marked points $p_1',\ldots,p_m'$, we can
insure that the corresponding evaluation map $\on{ev}^m: \ \M_d \to
G(r,n)^m$ is injective when restricted to $\sig_{I_*}(p_*).$ Let $Y
\subset G(r,n)^{l}\times G(r,n)^{m}$ be the closure of
$(\on{ev}^l \times \on{ev}^m)(\M_{d})$, and let $\phi: \
G(r,n)^{l}\times G(r,n)^{m} \to G(r,n)^m $ be the projection.

Since the homology class $[\phi(Y \cap \sigma_{I_*})]$ is non-trivial
\cite[page 64]{gr:pr}, $\phi(Y \cap \sigma_{I_*})$ must intersect
some Schubert variety
$$
\sigma_{J_*} = \sigma_{J_1} \times \sigma_{J_2} \times \ldots
\times \sigma_{J_m} \subset G(r,n)^m $$
of complementary dimension.  By Kleiman's lemma, \cite[Theorem 10.8
page 273]{ha:al}, the singular locus of $\phi(Y \cap \sigma_{I_*})$
does not intersect a general translate of $\sigma_{J_*}$, and
similarly the singular locus of $\sigma_{J_*}$ does not intersect
$\phi(Y \cap \sigma_{I_*})$.  Therefore the intersection occurs in the
smooth loci of $\phi(Y \cap \sigma_{I_*})$ and $\sigma_{J_*}$, and
another application of the lemma implies that the intersection is
finite.  

For generic translates of $\sigma_{J_*}$, the intersection is
contained in $\on{ev}^m(\sigma_{I_*}(p_*))$.  Indeed, let $
{\ol{\sig}}_{I_*}(p_*) $ be the compactification of $\sig_{I_*}(p_*)$
given in \cite{be:qs}, and $\Gamma \subset \ol{ \sig}_{I_*}(p_*)
\times G(r,n)^m$ the closure of the graph of $\on{ev}^m$.  Let $Z
\subset \Gamma$ be the complement of the graph of $\on{ev}^m$.  The
projection $\pi(Z)$ of $Z$ in $ G(r,n)^m$ is a closed sub-variety of
$\phi(Y \cap \sigma_{I_*})$.  By Kleiman's lemma, for generic
translates of $\sig_{J_*}$ the intersection of $\pi(Z)$ and
$\sig_{J_*}$ is empty, so the intersection is contained in
$\on{ev}^m(\sigma_{I_*}(p_*))$.

Because $\on{ev}^m \ | \ \sigma_{I_*}(p_*)$ is injective, the
intersection $\sig_{I_*}(p_*) \cap \sig_{J_*}(p_*') $ is finite and
non-empty.  Since the homology class $[\phi(Y \cap \sigma_{I_*})]$ is
independent of the choice of general translate of $\sigma_{I_*}$, the
above intersection is finite and non-empty for general translates of
the $\sig_{I_i}$ and $\sig_{J_j}$.  This implies that Gromov-Witten
invariant
$$
\l
[\sigma_{I_1}],\ldots,[\sig_{I_l}],[\sig_{J_1}],\ldots,[\sig_{J_m}]
\r_d \neq 0 .$$
In terms of the quantum product
$$
[\sigma_{I_1}] \star \ldots \star [\sig_{I_l}] \star [\sig_{J_1}]
\star \ldots \star [\sig_{J_{m-1}}] \supset q^d \sig_{*J_m} $$
which implies that $ [\sigma_{I_1}] \star \ldots \star [\sig_{I_l}]$
contains a term with coefficient $q^i$ with $i \leq d$.   
That is, for some Schubert variety $\sigma$, 
$$
\l [\sig_{I_1}], [\sig_{I_2}], \ldots , [\sig_{I_l}],[\sigma] \r_i
\neq 0 .$$
To prove the lemma it suffices to show that $i=d$.
By \cite[Moving Lemma 2.2]{be:qs}, for general translates of the
Schubert varieties the degree $i$ moduli space $\sig_{I_1}(p_1) \cap
\ldots \cap \sig_{I_l}(p_l) \cap \sigma(p)$ is finite and consists of
maps sending $p_1,\ldots,p_l,p$ to the corresponding Schubert cells.
Since $d$ is minimal, $i = d$.
\end{proof}

\section{Factorization}

In this section we show that a relationship between the polytopes for
different numbers of marked points is related to factorization of
Gromov-Witten invariants (i.e. associativity of quantum
multiplication).  A similar, easier, discussion holds for the additive
polytopes $\Del(l)$.  A consideration of a ``trivial'' factorization
completes the proof of Theorem \ref{final}.

Suppose that $l$ can be written $l = j + k -2 $ for
positive integers $j,k \ge 2$.  It is easy to see that $\Del_q(l)$ are 
projections of a section of $\Del_q(j) \times \Del_q(k)$\footnote{In fact, the volume functions
satisfy the factorization properties
$$ \Vol(\N(\mu_1,\ldots,\mu_{j-1},\nu_1,\ldots,\nu_{k-1}))
= \int_{\t_+}  \Vol(\N(\mu_1,\ldots,\mu_{j-1},*\lam)) 
\Vol(\N(* \lam, \nu_1,\ldots,\nu_{k-1}) )\d \lam $$
$$ \Vol(\M(\mu_1,\ldots,\mu_{j-1},\nu_1,\ldots,\nu_{k-1}))
= \int_{\Alc}  \Vol(\M(\mu_1,\ldots,\mu_{j-1},*\lam) )
\Vol(\M(* \lam, \nu_1,\ldots,\nu_{k-1})) \d \lam .$$
The second formula is implicit in Witten \cite[p.51]{wi:tw}, proved in
 \cite{jw:va}, and generalized in \cite{me:lo}.}.
$$  \Del_q(l) = \{ (\mu_1,\ldots,\mu_{j-1},\nu_1,\ldots,\nu_{k-1}) \ | \ 
(\mu,\nu) \in \Del_q(j) \times \Del_q(k), \ \ \mu_j = *\nu_k \} 
$$
To show the forward inclusion, note that if $A_1 A_2 \ldots A_l = I$
then letting $ B = A_j A_{j+1} \ldots A_l$ we have
$$(\lam(A_1),\ldots,\lam(A_{j-1}),\lam(B)) \in
\Del_q(j), $$
$$ (\lam(B^{-1}),\lam(A_j), \ldots, \lam(A_l) ) \in \Del(k).$$

In particular this means that any face of $\Del_q(l)$ is a projection
of a face (usually not of codimension $1$) of $\Del_q(j) \times
\Del_q(k)$.  Any face is the intersection of codimension $1$ faces.
This shows that any defining inequality of $\Del_q(l)$ is implied by a
finite set of defining inequalities for $\Del_q(j)$ and $\Del_q(k)$.

Using associativity of quantum cohomology one can be more specific
about which inequalities for $\Del_q(j),\Del_q(k)$ are needed to imply
an inequality for $\Del_q(l)$.  Suppose that a Gromov-Witten invariant
$ \l \sigma_{I_1},\ldots,\sigma_{I_l}, \sigma_J \r_d \neq 0 $
so that one has an inequality for $\Del(l)$ given by
\begin{equation} \labell{desired}
\sum_{i \in I_1} \lam_{1,i} + \ldots + \sum_{i \in I_{l}}
\lam_{l,i} \leq d .\end{equation}
Associativity of quantum multiplication says that 
$$ \l \sigma_{I_1},\ldots,\sigma_{I_l}, \sigma_J \r_d 
= \sum_{d_1 + d_2 = d, \ \ |K|=r  }
 \l \sigma_{I_1},\ldots,\sigma_{I_{j-1}}, \sigma_K \r_{d_1} 
 \l \sigma_{*K}, \sigma_{I_j},\ldots,\sigma_{I_l}, \sigma_J \r_{d_2} .$$
In particular there exist some $d_1,d_2$ with $d_1 + d_2 = d$ and
some Schubert variety $\sigma_K$ such that
$$  \l \sigma_{I_1},\ldots,\sigma_{I_{j-1}}, \sigma_K \r_{d_1} \neq 0, \ \ \
 \l \sigma_{*K} , \sigma_{I_j},\ldots,\sigma_{I_l}, \sigma_J \r_{d_2} \neq 0 .$$
From the non-vanishing of these Gromov-Witten invariants one deduces the 
inequalities for $\Del_q(j),\Del_q(k):$
\begin{equation} \labell{add1}
  \sum_{i \in I_1} \mu_{1,i} + \ldots + \sum_{i \in I_{j-1}}
\mu_{j-1,i} + \sum_{k \in K} \mu_{j,k} \leq d_1 ;\end{equation}
\begin{equation} \labell{add2}
 \sum_{k \in *K} \nu_{1,k} + \sum_{i \in I_{j}} \nu_{2,i} +
\ldots + \sum_{i \in I_{l}} \nu_{k,i} \leq d_2 .\end{equation}
Restricting to the section $\mu_j = *\nu_1$ one has that
$$ \sum_{k \in *K} \nu_{1,k} = - \sum_{k \in K} (* \nu_1)_k =
- \sum_{k \in K} \mu_{j,k}, $$
so by adding the two inequalities one obtains \eqref{desired}.

Using the trivial factorization $l = (l + 2) - 2$ we complete the
proof of Theorem \ref{final}.
\begin{lemma} \labell{factor}
Any inequality for $\Del_q(l)$ corresponding to a Gromov-Witten
invariant
$$ \l [\sigma_{I_1}],...,[\sigma_{I_l}],[\sigma_K]\r_d \neq 0 $$
is a consequence of an inequality corresponding to a Gromov-Witten
invariant of the form
$$ \l [\sigma_{I_1}],...,[\sigma_{I_{l-1}}],[\sigma_J]\r_{d_{1}} \neq 0 $$
for some $J \subset \{ 1,\ldots, n\}$ and $d_1\leq d.$
\end{lemma}

\begin{proof} Suppose that 
  $$  \l [\sigma_{I_1}],...,[\sigma_{I_l}],[\sigma_K]\r_d \neq 0 .$$
Taking $k=2$ we obtain that for some $J$ and $d_1\leq d$
$$ \l [\sigma_{I_1}],...,[\sigma_{I_{l-1}}],[\sigma_J]\r_{d_{1}}
\l [\sigma_{*J}],[\sigma_{I_l}],[\sigma_K]\r_{d_2} \neq 0 .$$
Thus the inequality
\begin{equation} \labell{want}
\sum_{i \in I_1} \lam_{1,i} + \ldots + \sum_{i \in I_{l}}
\lam_{l,i} \leq d .\end{equation}
follows from the inequalities
\begin{equation} \labell{first}
\sum_{i \in I_1} \lam_{1,i} + \ldots + \sum_{i \in I_{l-1}}
\lam_{l-1,i} +  \sum _{j\in J} \lambda_{l,j} \leq d_1 \end{equation}
for $\lam \in \Del_q(l)$ and
\begin{equation} \labell{last}
 \sum _{j\in *J} 
(*\lam_l)_j + \sum _{i\in I_l} (\lam_l)_i \leq d_2.
\end{equation}
The last equation is a tautology for $\lam_l \in \Alc$ by the $l=2$
case of Theorem \ref{result}.  In other words, \eqref{last} is implied
by the equations $\lam_{l,i} \ge \lam_{l,i+1}, \ \lam_{l,1} -
\lam_{l,n} \le 1$.  Thus \eqref{want} follows from \eqref{first} and
the inequalities defining $\Alc^l$.
\end{proof}

\section{Hidden symmetry}
\labell{symmetry}

An interesting aspect of the multiplicative problem is that it
possesses a symmetry not present in the additive case, related to the
symmetry of the fundamental alcove $\Alc$ of $SU(n)$.  Let
$Z \cong \Z/n\Z$ denote the center of $SU(n)$, with generator $c \in SU(n)$ 
the unique element of $SU(n)$ with
$$ \lambda(c) = (1/n,1/n, \ldots, 1/n, (1-n)/n) .$$
The action of $Z$ on $SU(n)$ induces an action on $\Alc \cong SU(n) /
\Ad(SU(n))$, given by
$$ c \cdot (\lambda_1,\ldots,\lambda_n) = (\lambda_2 + 1/n, \lambda_3
+ 1/n,\ldots, \lambda_n + 1/n, \lambda_1 - (n-1)/n) .$$

Let $C(l) \subset SU(n)^l$ denote the subgroup
$$ C(l) = \{ (z_1,\ldots,z_l) \subset Z^l \ | \ z_1z_2 \ldots z_l = 1 \} 
\cong Z^{l-1}.$$
The action of $C(l)$ on $\Alc^l$ leaves the polytope $\Del_q(l)$
invariant.

This symmetry of the polytope $\Del_q(l)$ implies a symmetry on the
facets of $\Del_q(l)$.  Let $c$ act on subsets of $\{ 1,2,\ldots,n \}$
via the action of $(12\ldots n)^{-1} \in S_n$:
$$ c^m \{ i_1, \ldots, i_r \} = \{i_{s+1} - m, \ldots, i_r - m,
 i_1 - m + n, \ldots, i_s - m + n\} $$
where $s$ is the largest index for which $i_s - m \ge 1$

Suppose an $l+1$-tuple $(I_1,\ldots,I_l,d)$ defines a facet of
$\Del_q(l)$ via the inequality \eqref{mult_ineq}.  Under the action of
$ (c^{m_1},\ldots,c^{m_l}) \in C(l)$, \eqref{mult_ineq}
becomes the inequality corresponding to $(c^{m_1}I_1,\ldots,c^{m_l}I_l,d')$
where $d'$ is defined by
\begin{equation} \labell{4ac}
 \sum_{i=1}^l | c^{m_i}I_i | + nd' =  \sum_{i=1}^l | I_i |  + nd.
\end{equation}
\begin{example}  From the degree $0$ inequality $\lambda_n(A) + 
  \lambda_n(B) \le \lambda_n(AB)$ we obtain by the action of
  $(c^{-i},c^{-j},c^{i+j}), \ i + j \le n$ the degree $1$ inequality
  \eqref{quant_Weyl}.
\end{example}

Equation \eqref{4ac} defines a $C(l)$ action on the set of
$l+1$-tuples $(I_1,\ldots,I_l,d)$ defining facets of $\Del_q(l)$.  It
is an interesting fact that the Gromov-Witten invariants $\l
\sig_{I_1},\ldots,\sig_{I_l}\r_d$ are invariant under this action:
\begin{proposition} Let $(c^{m_1},\ldots,c^{m_l}) \in C(l)$.  Then
$\l \sig_{I_1},\ldots,\sig_{I_l} \r_d = \l \sigma_{c^{m_1}
I_1}, \ldots, \sigma_{c^{m_l} I_l} \r_{d'}.$
\end{proposition}
\begin{proof}  Let $\sigma_c = \sigma_{r,r+1,\ldots,n-1}$ denote the 
Schubert variety isomorphic to the Grassmannian $G(r,n-1)$ of
$r$-planes contained in $n-1$-space.  We claim that quantum
multiplication by $\sigma_c$ is given by the following formula:
\begin{equation}  \labell{Cox_mult}
[ \sigma_c ] \star [ \sigma_I ] = q^{(|cI| + r - |I|)/n} [ \sigma_{cI} ] .
\end{equation}
The exponent $(|cI| + r - |I|)/n$ equals $1$ if $1 \in I$, and equals
$0$ otherwise. In particular $[\sig_c]^{\star n}=q^r.$ 

The lemma then follows by associativity of the quantum product.
Without loss of generality it suffices to show that the Gromov-Witten
invariants are invariant under an element of the form
$(c,c^{-1},1,\ldots,1) \in C(l)$.  
Given that 
$$
[\sig _{I_1}]\star \ldots \star [\sig _{I_{l-1}}]\supset
\l \sig_{I_1},\ldots,\sig_{I_l} \r_d
[\sigma_{*I_{l}}]q^d$$
multiplying by $[\sig_c]$ on both sides yields
$$[\sig _{cI_1}]\star \ldots \star [\sig _{I_{l-1}}]\supset \l
\sig_{I_1},\ldots,\sig_{I_l} \r_d [\sigma_{c(*I_{l})}]q^{d'}=\l
\sig_{I_1},\ldots,\sig_{I_l} \r_d [\sig_{*c^{-1}I_l}]q^{d'}.
$$

The formula \eqref{Cox_mult} may be proved using either the canonical
isomorphism of quantum Schubert calculus with the Verlinde algebra of
$U(r)$,
$$ QH^*(G(r,n))/(q=1) \cong R(U(r)_{n-r,n}).$$
given a mathematical proof in Agnihotri \cite{ag:th}, or
using the combinatorial formula of Bertram, Ciocan-Fontanine and
Fulton \cite{be:qm}.  $R(U(r)_{n-r,n})$ denotes the Verlinde
algebra of $U(r)$ at $SU(r)$ level $n-r$ and $U(1)$ level $n$, and is
the quotient of the tensor algebra $R(U(r))$ by the relations
$$ V_{\lam} \sim (-1)^l(w) V_{w(\lam+ \rho)-\rho}, w \in \Waff $$
and if $\lam_1 - \lam_r \le n -r $ then
$$ V_{(\lam_1,\ldots,\lam_r)} \sim V_{(\lam_2-1,\lam_3-1,
\ldots,\lam_r-1,\lam_1 - (n-r+1)} .$$
Here $\Waff$ acts on $\Lambda^*$ at level $n$, and $\rho$ is the
half-sum of positive roots.  The Verlinde algebra $R(U(r)_{n-r,n})$
has as a basis the (equivalence classes of the) representations
$V_{\lam}$, where $\lam = (\lambda_1,\ldots,\lambda_r) \in \Z^r, \ 0
\leq \lambda_i \leq n- r $ are dominant weights of $U(r)$ at level
$n-r$.

The canonical isomorphism is given by $\sigma_I \mapsto V_{\lambda}$,
where $\lambda$ is defined by
$$ \lambda_j = n-r + j - i_j .$$

The key point is that the sub-algebra $R(U(1)) \subset R(U(r))$
descends to a sub-algebra $ R(U(1)_n) \subset R(U(r)_{n-r,r})$
generated by the representation $V_c := V_{(1,1,\ldots,1)}$, which 
maps under the isomorphism to the Schubert variety $\sigma_c$.
From the description of the algebra given above one sees
that 
$ V_c \fus V_{\lam} = V_{\lambda'} $
where
$$ \lambda' = 
\begin{array}{cl}
(\lambda_1 + 1,\lambda_2 +1,\ldots, \lambda_r + 1) &  \hbox{\ if\ }\lambda_1 < n-r
\\ (\lambda_2,\ldots,\lambda_r,\lambda_1 - n +r) &  \hbox{\ if\ }\lambda_1 = n-r
\end{array} .$$
Since $V_{\lambda'}$ maps to $\sigma_{cI}$ under the canonical
isomorphism, this proves \eqref{Cox_mult}.

Alternatively, \eqref{Cox_mult} can be derived from the combinatorial
rim-hook formula of \cite[p. 8]{be:qm}.  Let $\lambda^t$ denote the
transpose of $\lambda$, so that $\sig_{\lambda^t}$ is the image of
$\sig_\lambda$ under the isomorphism $G(r,n) \cong G(n-r,n)$.  The
ordinary (resp. quantum) Littlewood-Richardson numbers are invariant
under transpose
$$ N_{\lam^t \mu^t}^{\rho^t} =  N_{\lam \mu}^{\rho}, \ \ \
 N_{\lam^t \mu^t}^{\rho^t}(n-r,r) =  N_{\lam \mu}^{\rho}(r,n-r).
 $$
It follows from \cite[Corollary]{be:qm} that 
$$   N_{\lam \mu}^{\rho}(r,n-r) =\sum \eps(\rho^t/\nu^t) 
N_{\lam \mu}^{\rho} $$
where $\rho$ ranges over all diagrams of height $\leq r$ that can be
obtained by adding $m$ rim-hooks.

If $\mu = (1,1,\ldots,1)$ then 
$$ V_\mu \otimes V_{\lam} = V_\rho, \ \ \rho =
(\lam_1 + 1, \ldots, \lam_r + 1) .$$
If $\lam_1 < n-r$, then since the height of $\rho$ is $\leq r$, there
are no rim $n$-hooks in $\rho$.  On the other hand, if $\lam_1 = n-r$,
then it is easy to see that there is a unique rim $n$-hook in $\rho$,
whose complement is $\lam'$ above.

We have learned from A. Postnikov that formula similar to
\eqref{Cox_mult} holds for the full flag variety \cite{po:hs}.  A
deeper reason for the appearance of symmetry is given by Seidel
\cite{se:pi}.
\end{proof}

This symmetry simplifies the computation of many Gromov-Witten
invariants.  For sufficiently small $n$ and $l$ all Gromov-Witten
invariants are equivalent to degree $0$ ones.  An example of a
Gromov-Witten invariant not equivalent via symmetry to a degree $0$
invariant is the degree $1$ invariant for $G(5,10)$
$$ \l \sigma_{\{ 2,4,6,8,10 \}}, \sigma_{\{ 2,4,6,8,10 \}}, \sigma_{\{
 1,3,5,7,9 \}} \r_1 $$
which may be computed using the formula of \cite{be:qm}.  That is, for
sufficiently large $r,n$, not all of the inequalities are related to
``classical'' inequalities via symmetry.

\section{Verlinde algebras}
\labell{Verlinde}

Finally we want to explain the representation-theoretic interpretation
of $\Del_q(l)$ in terms of the Verlinde algebra of $SU(n)$.  Denote by
$\Lambda^*_N$ the set of dominant weights of $SU(n)$ at level $N$:
$$ \Lambda^*_N = \{ (\lambda_1,\ldots,\lambda_n) \in (\Z/n)^n \ | \
\lam_i - \lam_{i+1} \in \Z_{\ge 0},\ \  \lam_1 - \lam_n \le N \} .$$
The Verlinde algebra $R(SU(n)_N)$ is the free group $\Z[\Lambda_N^*]$
on the generators $V_{\xi}, \xi \in \Lambda_N^*$.  The algebra
structure is given by ``fusion product''
$$ V_{\xi_1} \fus_N \ldots \fus_N \V_{\xi_l} = \sum_{\nu \in
\Lambda^*_N}  m^N({\xi_1,\ldots,\xi_l,\nu}) \  V_{* \nu} $$
where the coefficients $m^N(\xi_1,\ldots,\xi_l)$ are defined as
follows.  There is a positive line bundle $L^N({\xi_1,\ldots,\xi_l})
\to\M(\xi_1/N ,\ldots ,\xi_l/N)$ which descends from the polarizing
line bundle on $\ti{R}$ (see Pauly \cite[Section 3]{pa:em}).  The
coefficient $m^N({\xi_1,\ldots,\xi_l})$ is defined by
$$ m^N({\xi_1,\ldots,\xi_l}) = \dim(H^0(L^N(\xi_1,\ldots,\xi_l)) .$$ 
Since $L^N$ is positive, $\M(\xi_1/N,\ldots,\xi_l/N)$ is non-empty if
and only if for some $k$
$$\dim(H^0(L^k(\xi_1,\ldots,\xi_l)^{\otimes
N})=\dim(H^0(L^{kN}(k\xi_1,\ldots,k\xi_l))\neq 0,$$ 
that is, $m^{kN}({k\xi_1,\ldots,k\xi_l}) \neq 0.$


\end{document}